\documentclass[a4paper,12pt]{amsart}
\usepackage{amssymb}
\usepackage{ifthen}
\nonstopmode \numberwithin{equation}{section}
\newtheorem{thm}{Theorem}
\newtheorem{cor}{Corollary}
\newtheorem{lem}{Lemma}

\newtheorem{conj}{Conjecture}

\theoremstyle{definition}
\newtheorem{defn}{Definition}[section]

\newtheorem{prob}[equation]{Problem}

\newenvironment{rem}{%
\bigskip
\noindent \textsl{{\sl Remark. }}}{\bigskip}
\newenvironment{rems}{%
\bigskip
\noindent \textsl{{\sl Remarks. }}}{\bigskip}

\newcounter {own}
\def\theown {\thesection       .\arabic{own}}

\newenvironment{pf}[1][]{%
 \vskip 3mm
 \noindent
 \ifthenelse{\equal{#1}{}}%
  {{\slshape Proof. }}%
  {{\slshape #1.} }%
 }%
{\qed\bigskip}

\newcounter{alphabet}
\newcounter{tmp}
\newenvironment{Thm}[1][]{\refstepcounter{alphabet}%
\bigskip%
\noindent%
{\bf Theorem \Alph{alphabet}}%
\ifthenelse{\equal{#1}{}}{}{ (#1)}%
{\bf .} \itshape}{\vskip 8pt}

\newcommand{\ID}{{\mathbb D}}

\newcommand{\IC}{{\mathbb C}}

\def\be{\begin{equation}}
\def\ee{\end{equation}}

\newcommand{\bee}{\begin{enumerate}}
\newcommand{\eee}{\end{enumerate}}

\newcommand{\blem}{\begin{lem}}
\newcommand{\elem}{\end{lem}}
\newcommand{\bthm}{\begin{thm}}
\newcommand{\ethm}{\end{thm}}
\newcommand{\bcor}{\begin{cor}}
\newcommand{\ecor}{\end{cor}}
\newcommand{\beg}{\begin{examp}}
\newcommand{\eeg}{\end{examp}}
\newcommand{\begs}{\begin{examples}}
\newcommand{\eegs}{\end{examples}}
\newcommand{\bdefe}{\begin{defn}}
\newcommand{\edefe}{\end{defn}}
\newcommand{\bprob}{\begin{prob}}
\newcommand{\eprob}{\end{prob}}
\newcommand{\bei}{\begin{itemize}}
\newcommand{\eei}{\end{itemize}}

\newcommand{\bcon}{\begin{conj}}
\newcommand{\econ}{\end{conj}}
\newcommand{\bcons}{\begin{conjs}}
\newcommand{\econs}{\end{conjs}}
\newcommand{\bprop}{\begin{propo}}
\newcommand{\eprop}{\end{propo}}
\newcommand{\br}{\begin{rem}}
\newcommand{\er}{\end{rem}}
\newcommand{\brs}{\begin{rems}}
\newcommand{\ers}{\end{rems}}
\newcommand{\bo}{\begin{obser}}
\newcommand{\eo}{\end{obser}}
\newcommand{\bos}{\begin{obsers}}
\newcommand{\eos}{\end{obsers}}
\newcommand{\bpf}{\begin{pf}}
\newcommand{\epf}{\end{pf}}
\newcommand{\ba}{\begin{array}}
\newcommand{\ea}{\end{array}}
\newcommand{\beq}{\begin{eqnarray}}
\newcommand{\beqq}{\begin{eqnarray*}}
\newcommand{\eeq}{\end{eqnarray}}
\newcommand{\eeqq}{\end{eqnarray*}}

\newcommand{\ra}{\rightarrow}

\newcounter{minutes}
\divide\time by 60
\newcounter{hours}
\multiply\time by 60 \addtocounter{minutes}{-\time}
\begin{document}
\title{Sufficient conditions for univalence and study of a class of meromorphic univalent functions}
\begin{center}
{\tiny \texttt{FILE:~\jobname .tex,
        printed: \number\year-\number\month-\number\day,
        \thehours.\ifnum\theminutes<10{0}\fi\theminutes}
}
\end{center}

\author{Bappaditya Bhowmik${}^{~\mathbf{*}}$}
\address{Bappaditya Bhowmik, Department of Mathematics,
Indian Institute of Technology Kharagpur, Kharagpur - 721302, India.}
\email{bappaditya@maths.iitkgp.ernet.in}
\author{Firdoshi Parveen}
\address{Firdoshi Parveen, Department of Mathematics,
Indian Institute of Technology Kharagpur, Kharagpur - 721302, India.}
\email{frd.par@maths.iitkgp.ernet.in}

\begin{abstract}
In this article we consider the class $\mathcal{A}(p)$ which consists of functions that
are meromorphic in the unit disc $\ID$ having a simple pole at $z=p\in (0,1)$ with the normalization $f(0)=0=f'(0)-1 $.
First we prove some sufficient conditions for univalence of such functions in $\ID$. One of these conditions enable us to
consider the class  $\mathcal{V}_{p}(\lambda)$ that consists of functions satisfying certain differential inequality which forces univalence of
such functions. Next
we establish that $\mathcal{U}_{p}(\lambda)\subsetneq \mathcal{V}_{p}(\lambda)$, where $\mathcal{U}_{p}(\lambda)$ was introduced
and studied in \cite{BF-1}. Finally, we discuss some coefficient problems for  $\mathcal{V}_{p}(\lambda)$ and end the article with a coefficient conjecture.
\end{abstract}
\thanks{}
\maketitle
\pagestyle{myheadings}
\markboth{B. Bhowmik and F. Parveen}{Sufficient conditions for univalence and study of a class of meromorphic univalent functions}

\bigskip
\section{Introduction and sufficient condition for univalence}

Let $\mathcal{M}$ be the set of meromorphic functions $F$ in $\Delta=\{\zeta\in \IC:|\zeta|>1\}\cup\{\infty\}$ with the following expansion:
$$
F(\zeta)=\zeta+\sum_{n=0}^{\infty} b_{n}\zeta^{-n},\quad \zeta\in\Delta.
$$
This means that these functions have simple pole at $z=\infty$ with residue $1$.
Let $\mathcal{A}$ be the collection of all analytic functions in $\ID$ with the normalization $f(0)=0=f'(0)-1 $.
In \cite{aksen}, Aksent\'{e}v proved a sufficient condition for a function $F\in \mathcal{M}$ to be univalent which we state now:

\begin{Thm}\label{TheoA}
If $F\in \mathcal{M}$ satisfies the inequality
$$
|F'(\zeta)-1|\leq1,\quad \zeta\in\Delta,
$$
then $F$ is univalent in $\Delta$.
\end{Thm}
This result motivated many authors to consider the classes
$\mathcal{U}(\lambda):=\{f\in \mathcal{A}: |U_f(z)|< \lambda\}, \lambda\in
(0,1]$ where $U_f(z):=(z/f(z))^2 f'(z)-1$ and this class has been studied
extensively in \cite{OP,OPW} and references therein. In \cite{BF-1}, we
wanted to see the meromorphic analogue of the class $\mathcal{U}(\lambda)$ by
introducing a nonzero simple pole for such functions in $\ID$. More
precisely, we consider the class $\mathcal{A}(p)$ of all functions $f$ that
are holomorphic in $\ID\setminus \{p\}$, $p\in(0,1)$  possessing a simple
pole at the point $z=p$ with nonzero residue $m$ and normalized by the
condition $f(0)=0=f'(0)-1 $. We define $\Sigma(p):=\{f\in\mathcal{A}(p): f
\mbox{ is one to one in} ~~\ID \}$. Therefore, each $f\in \mathcal{A}(p)$ has
the Laurent series expansion of the following form
\beq \label{fp4eq1}
f(z)=\frac{m}{z-p}+\sum\limits_{n=0}^\infty a_n z^n,~~z\in \ID\setminus\{p\}.
\eeq
In this context we proved a sufficient condition for a function $f\in
\mathcal{A}(p)$ to be univalent (see \cite[Theorem 1]{BF-1}), which we recall
now.

\begin{Thm}\label{TheoB}
Let $f\in \mathcal{A}(p)$. If $\left|U_f(z)\right|\leq\left((1-p)/(1+p)\right)^2$ for $z\in \ID$,
then $f$ is univalent in $\ID$.
\end{Thm}
Using Theorem~B, we constructed a subclass $\mathcal{U}_{p}(\lambda)$ of $\Sigma(p)$ which is defined as follows:
$$
\mathcal{U}_p(\lambda):=\left\{f\in \mathcal{A}(p):\left|U_{f}(z)\right|<\lambda \mu, z\in \ID\right\}
$$
where $0<\lambda \leq 1$ and $\mu=((1-p)/(1+p))^2$. In this note, we improve
the sufficient condition proved in Theorem~B by replacing the number $\mu=((1-p)/(1+p))^2$ with the
number $1$. We give a proof of this result below.
\bthm\label{fp4th1}
Let $f\in \mathcal{A}(p)$. If $\left|U_{f}(z)\right|<1$
holds for all $z\in \ID$ then $f\in \Sigma(p)$.
\ethm
\bpf
Let $\mathcal{M}_p:=\{ f\in \mathcal{M}: F(1/p)=0\}$ where $0<p<1$. Clearly,
$\mathcal{M}_p \subseteq \mathcal{M}$. For each $f\in\mathcal{A}(p)$ consider the
transformation $F(\zeta):=1/f(1/\zeta)$, $\zeta \in \Delta$. We claim that $F\in\mathcal{M}_p \subseteq \mathcal{M}$. Since $f$ has an
expansion of the form (\ref{fp4eq1}), therefore we have
\beqq
F(\zeta)&=&1/f(1/\zeta)\\
&=& \left(m \zeta/(1-p\zeta)+\sum_{n=0}^{\infty}a_n\zeta^{-n}\right)^{-1}\\
&=&\zeta+ (a_1-p a_2-1)/p +\left(p(a_2-pa_3)+(a_1-pa_2)^2-(a_1-pa_2)\right)/\zeta p^2+\cdots .
\eeqq
Here we see that
$F(1/p)=0, F(\infty)=\infty$ and $F'(\infty)=1$.
This proves that each $f\in \mathcal{A}(p)$ can be associated with the mapping $F\in\mathcal{M}_p $.
 Using the change of variable $\ID \ni z=1/\zeta$, the above association quickly yields
$$
F'(\zeta)-1=f'(1/\zeta)/(\zeta^{2}f^{2}(1/\zeta))-1=z^{2}f'(z)/f^{2}(z)-1=U_{f}(z).
$$
Now since $\mathcal{M}_p \subseteq \mathcal{M}$, an application of the Theorem~A gives that if any function $F\in\mathcal{M}_p $ satisfies
$|F'(\zeta)-1|\leq1,\quad \zeta\in\Delta$, then $F$ is univalent in $\Delta$, i.e. the inequality $|U_{f}(z)|<1$ forces $f$ to be univalent in $\ID$.
\epf

In view of the Theorem~\ref{fp4th1}, it is natural to consider a new subclass $\mathcal{V}_{p}(\lambda)$ of $\Sigma(p)$ defined as:
$$
\mathcal{V}_p(\lambda):=\left\{f\in \mathcal{A}(p):\left|U_{f}(z)\right|<\lambda, ~z\in \ID\right\},\quad \mbox{for}~~ \lambda\in(0,1].
$$
We now claim that $\mathcal{U}_{p}(\lambda)\subsetneq \mathcal{V}_{p}(\lambda) \subsetneq \Sigma(p)$. To establish the first inclusion,
we note that as $\lambda \mu < \lambda$, therefore we have $\mathcal{U}_{p}(\lambda)\subseteq \mathcal{V}_{p}(\lambda)$. Now consider the function
$$
k_p^{\lambda}(z):=\frac{-pz}{(z-p)(1-\lambda pz)}, z\in \ID.
$$
It is easy to check that $ U_{k_p^{\lambda}}(z)=-\lambda z^2$ so that $| U_{k_p^{\lambda}}(z)|<\lambda$ but
$|U_{k_p^{\lambda}}(z)|\nless \lambda\mu$ for all $z\in \ID$. This proves the first inclusion. Next we wish to establish the second inclusion of our claim.  We see that by virtue of
the Theorem~\ref{fp4th1}, $\mathcal{V}_{p}(\lambda)\subseteq \Sigma(p)$. Again considering the following two examples,
we see that $\mathcal{V}_{p}(\lambda)\subsetneq \Sigma(p)$ for $0<\lambda\leq 1$.\\
\textbf{Case 1:}\,($0<\lambda<1$). Take $a\in \IC$ such that $\lambda<|a|<1$. Consider the functions $f_a$ defined by
 $$
f_a(z)=\frac{z}{(z-p)(az-1/p)}, \quad z\in \ID.
$$
 It is easy to check that $f_a$ satisfies the normalizations $f_a(p)=\infty$ and $f_a(0)=0=f_a'(0)-1$. Also $f_a(z)$ is univalent in $\ID$ and  $U_{f_a}(z)=-az^2$. Now as $|z|\ra 1^{-}$, $|U_{f_a}(z)|\ra |a| > \lambda$. Therefore $f_a(z)\notin \mathcal{V}_{p}(\lambda)$. This shows that $\mathcal{V}_{p}(\lambda)$ is a proper subclass of $\Sigma(p)$ for $0<\lambda<1$. \\
 \textbf{Case 2:}\,($\lambda=1$). It is well-known that the function
$$
g(z)=\frac{z-\frac{2p}{1+p^2}z^2}{(1-\frac{z}{p})(1-zp)}, z\in \ID,
$$
is in $\Sigma(p)$ (Compare \cite{BPW}). A little calculation shows that
$$
U_g(z)=\left(z(1-p^2)/(1+p^2)\right)^2\left(1-(2pz/(1+p^2))\right)^{-2}.
$$
Now $|U_g(z)|<1$ holds for all $|z|\leq R$ whenever $R< \frac{1+p^2}{1+2p-p^2}<1$. From here we can conclude that $g$ does not belongs to the class $\mathcal{V}_{p}(\lambda)$ for $\lambda=1$, i.e. $\mathcal{V}_{p}:=\mathcal{V}_{p}(1) \subsetneq \Sigma(p)$.

\br
It can be easily seen that similar to the class $\mathcal{U}_{p}(\lambda)$, the class $\mathcal{V}_{p}(\lambda)$ is preserved under conjugation
and is not preserved under the operations like rotation, dilation, omitted value transformation and the $n$-th root transformations.
\er

Let $f\in\mathcal{A}(p)$. We see that the function $z/f$ is analytic in $\ID$ and non vanishing in $\ID\setminus\{p\}$. Therefore it has a Taylor
expansion of the following form about the origin.
\beq\label{fp4eq5}
 \frac{z}{f(z)}=1+b_{1}z+b_{2}z^2+\cdots, ~~~z\in \ID.
\eeq
Now we prove some sufficient conditions for univalence of functions  $f\in\mathcal{A}(p)$ which involves the second and higher order derivatives of $z/f$.
These are the contents of the next two theorems.
\bthm\label{fp4th3}
Let  $f\in\mathcal{A}(p)$ and $f/z$ is non-vanishing in $\ID \setminus\{0\}$. If $|(z/f(z))''|\leq 2$ for $z\in \ID$, then $f$ is univalent in $\ID$.
This condition is only sufficient for univalence but not necessary.
\ethm
\bpf
First we prove the univalence of $f$. Using the expansion (\ref{fp4eq5}), we have
$$
U_f(z)=-z(z/f)'+(z/f)-1=\sum_{n=2}^{\infty}(1-n)b_n z^n.
$$
We also note that $zU_f'(z)=-z^2(z/f)''$. Therefore $|(z/f)''|\leq 2$  yields $|zU_f'(z)|\leq 2|z|$. This implies that $zU_f'(z)\prec 2z$ where $\prec$ denotes usual
subordination. Now by a well known result of subordination (compare \cite[p. 76, Theorem 3.1d.]{MM}), we get $U_f(z)\prec z$, i.e. $|U_f(z)|\leq |z|<1$.
This shows that $f$ is univalent in $\ID$ by virtue of the Theorem~\ref{fp4th1}.
In order to establish the second claim of the theorem, we consider the function
$$
h(z)= \frac {2pz}{(p-z)(2-pz(p+z))},\, z\in \ID.
$$
Note that $h(0)=0=h'(0)-1$ and $h(p)=\infty$. Also since $|pz(p+z)|<2$,  $h$ has no other poles in $\ID$ except at $z=p$.
Consequently $h\in \mathcal{A}(p)$.
It is easy to check that $U_{h}(z)=-z^3$ and $(z/h)''= 3 z$. Hence $|U_{f}(h)|<1$  but $|(z/h)''|> 2$ for $2/3<|z|<1$.
This example shows that the boundedness condition in the statement of the Theorem is only sufficient but not necessary.
\epf

The following theorem is also a univalence criteria described by a sharp inequality involving the
$n$-th order derivatives  of $z/f$ (denoted by $(z/f)^n$), $n\geq 3$.
\bthm\label{fp4th4}
Let $f\in \mathcal{A}(p)$ and $f(z)\neq 0$ for $\ID\setminus\{0\}$. If for $n\geq 3$,
\beq\label{fp4eq8}
\sum_{k=0}^{n-3}\frac{k+1}{(k+2)!}|\alpha_k|+\frac{n-1}{n!}\left|\left(\frac{z}{f}\right)^n\right|\leq 1,~ z\in \ID,
\eeq
where $\alpha_k=- (z/f)^{k+2}|_{z=0}$, then $f$ is univalent in $\ID$. The result is sharp and equality holds in the above inequality for the function
$k_p(z)=-pz/(z-p)(1-pz)$ for all $n\geq 3$ and for the functions
$$
f_n(z)= \frac {z}{1-\left(1/p+ p^{n-1}/(n-1)\right)z+ z^n/(n-1)}, \, z\in \ID,
$$
for each $n\geq 3$.
\ethm
\bpf
Proceeding similarly as  the proof of \cite[Theorem 1.1]{OP1}, the inequality (\ref{fp4eq8}) will imply that $|U_f(z)|<1$
which proves that $f$ is univalent in $\ID$.
To complete the proof of remaining assertion of the theorem, we consider the univalent function $k_p(z)$, $z\in \ID$ and compute
$$
(z/k_p(z))'=-(1/p+p)+2z,\quad (z/k_p(z))''=2 \quad \mbox{and} \quad (z/k_p(z))^n=0,\,n\geq3.
$$
Therefore we get $\alpha_0=-2$ and $\alpha_k=0$ for $k\geq 1$.
Taking account of the above computations, it can now be easily checked that the equality holds in the inequality (\ref{fp4eq8}).
Lastly, It can be proved that the functions $f_n \in \mathcal{V}_{p}(\lambda)$ for $\lambda=1$ i.e., $f_n$ is univalent in $\ID$.
Again for the  functions $f_n$, it is easy to check that
$\alpha_k=0,\,0\leq k\leq n-3$   and $(z/f_n)^n=n!/(n-1)$ for all $n\geq 3$, which essentially proves the sharpness of the result.
\epf

Now in the following theorem we give sufficient conditions for a function $f\in \mathcal{A}(p)$ to be in the class $\mathcal{V}_{p}(\lambda)$
by using Theorem~\ref{fp4th1}, Theorem~\ref{fp4th3} and Theorem~\ref{fp4th4} in terms of the coefficients $b_n$ defined in (\ref{fp4eq5}).
\bthm\label{fp4th5}
Let $f\in \mathcal{A}(p)$ and each $z/f$ has the expansion of the form $(\ref{fp4eq5})$. If $f$ satisfies any one of the following three conditions namely \\
$(i) ~~\sum_{n=2}^{\infty}(n-1)|b_n|\leq \lambda$\\
$(ii)~~ \sum_{n=2}^{\infty}n(n-1)|b_n|\leq 2\lambda$\\
$(iii)~~ \sum_{k=2}^{n}(k-1)|b_k|+(n-1)\sum_{k=n+1}^{\infty}\binom{k}{n}|b_k|\leq \lambda$\\
then $f\in \mathcal{V}_{p}(\lambda)$.
\ethm
\bpf
Since $z/f$ has the form (\ref{fp4eq5}), it is simple exercise to see that
$$
U_f(z)=-\sum_{n=2}^{\infty}(n-1)b_n z^n,\quad (z/f)''=\sum_{n=2}^{\infty}n(n-1)b_n z^{n-2}
$$
and
$$
\left(\frac{z}{f}\right)^n= n! b_n+\sum_{k=n+1}^{\infty}\frac{k!b_k}{(k-n)!} z^{k-n}=\sum_{k=n}^{\infty}\frac{k!b_k}{(k-n)!} z^{k-n}
$$
Therefore condition (i) and (ii) implies that $|U_f(z)|<\lambda$ and $|(z/f)''|<2 \lambda$ respectively.
Again following the similar arguments of the proof of the Theorem~\ref{fp4th3}, we conclude that $|(z/f)''|<2 \lambda$ implies $|U_f(z)|<\lambda$.
Now
$$\alpha_k=- (z/f)^{k+2}|_{z=0}=-(k+2)! b_{k+2}.
$$
Substituting the value of $\alpha_k$ and $(z/f)^n$ in terms of the coefficient $b_n$ in the left hand side of the inequality (\ref{fp4eq8}) we get
\beqq
&& \sum_{k=0}^{n-3}(k+1)|b_{k+2}|+\frac{n-1}{n!}\left|\sum_{k=n}^{\infty}\frac{k!b_k}{(k-n)!} z^{k-n}\right|\\
&\leq& \sum_{k=2}^{n}(k-1)|b_k|+(n-1)\sum_{k=n+1}^{\infty}\binom{k}{n}|b_k| \\
&\leq& \lambda \quad \mbox{(by (iii))}
\eeqq
Hence an application of the Theorem~\ref{fp4th4} gives $|U_f(z)|<\lambda$.
This shows that in each case $f\in \mathcal{V}_{p}(\lambda)$.
\epf

In the following section we study some coefficient problem for functions in $\mathcal{V}_{p}(\lambda)$ which is one of the important problem in geometric function theory.

\section{Coefficient problem for the class $\mathcal{V}_{p}(\lambda)$}

Let $f\in\mathcal{V}_{p}(\lambda)$ with the expansion (\ref{fp4eq5}). Now proceeding as a similar manner of ( \cite[Theorem 12]{BF-1}) we have the sharp bounds for $|b_n|, n\geq 2$, which is given by
$$
|b_n|\leq \frac{\lambda }{n-1}, \quad n\geq 2,
$$
and equality holds in the above inequality for the function
\beq\label{fp4eq6}
f(z)= \frac {z}{1-\left(1/p+(\lambda p^{n-1})/(n-1)\right)z+\lambda z^n/(n-1)}, \, z\in \ID.
\eeq

Each $f\in \mathcal{V}_{p}(\lambda) $ has the following Taylor expansion
\beq\label{fp4eq7}
f(z)=z+\sum_{n=2}^{\infty}a_{n}(f)z^{n},\quad |z|<p.
\eeq

Now the problem is to find out the region of variability of these Taylor coefficients $a_n, n\geq2 $. Here we note that similar to the class $\mathcal{U}_{p}(\lambda)$, every $f\in \mathcal{V}_{p}(\lambda)$ has the following representation (see \cite[Theorem 3]{BF-1}):
\beq\label{fp4eq3}
\frac{z}{f(z)}=1-\left(\frac{f''(0)}{2}\right)z+\lambda z\int_{0}^{z}w(t)dt,
\eeq
where $w\in\mathcal{B}$. Here $\mathcal{B}$ denotes the class of functions $w$ that are analytic in $\ID$ such that $|w(z)|\leq1$ for $z\in \ID$. By using this representation formula in the following theorem we give the exact set of variability for the second Taylor coefficient of $f\in\mathcal{V}_{p}(\lambda)$.

\bthm\label{fp2th2}
Let each $f\in\mathcal{V}_{p}(\lambda)$ has the Taylor expansion $f(z)=z+\sum_{n=2}^{\infty}a_{n}(f)z^{n},$ in the disc $\{z: |z|<p\}.$
Then the exact region of variability of the second Taylor coefficient $a_{2}(f)$ is the disc determined by the inequality
\beq
\label{fp4eq4}
|a_{2}(f)-1/p| \leq \lambda p.
\eeq
\ethm

\bpf
Substituting $z=p$ in (\ref{fp4eq3}) we get
$$
a_2(f)=\frac{f''(0)}{2}=\frac{1+\lambda p\int_{0}^{p}w(t)dt}{p}
$$
which implies
\beqq
|a_{2}(f)-1/p|&=& \left|\frac{\lambda p\int_{0}^{p}w(t)dt}{p}\right|\\
&\leq& \lambda \int_{0}^{p}|w(t)|dt \leq \lambda p.
\eeqq
Therefore $|a_{2}(f)-1/p|\leq \lambda p$.
A point on the boundary of the disc described by (\ref{fp4eq4}) is attained for the function
$$
f_{\theta}(z)=\frac{z}{1-\frac{z}{p}\left(1+\lambda p^{2}e^{i\theta}\right)+\lambda e^{i\theta}z^2}
$$
where $\theta\in [0,2\pi].$ Also the points in the interior of the disc described in (\ref{fp4eq4}) are attained by the functions
$$
f_{a}(z)=\frac{z}{1-\frac{z}{p}\left(1+\lambda a p^{2}\right)+\lambda a z^2}
$$
where $0<|a|<1$. It is easy to see that these functions belong to the class $\mathcal{V}_{p}(\lambda)$.
This shows that the exact region of variability of $a_{2}(f)$ is given by the disc (\ref{fp4eq4}).
\epf

Following consequences of the above theorem can be observed easily:
\bcor
Let for some $\lambda \in (0,1]$, $f\in\mathcal{V}_{p}(\lambda)$ and has the form $f(z)=z+\sum_{n=2}^{\infty}a_{n}(f)z^{n}$, in $|z|<p$.
Then $|a_{2}(f)|\leq 1/p+\lambda p$ and equality holds in this inequality for the function
$k_p^{\lambda}$.
\ecor
Now the function $k_p^{\lambda}$ is analytic in the disk $\{z: |z|<p\}$ and has the Taylor expansion as
$$
k_p^{\lambda}(z)=\sum_{n=1}^{\infty}\frac{1-\lambda^n p^{2n}}{p^{n-1}(1-\lambda p^2)}z^n, \quad |z|<p.
$$
Since the function $k_p^{\lambda}$ serves as an extremal function for the class  $\mathcal{V}_{p}(\lambda)$,  the above corollary enables us to make the following
\bcon
If $f\in\mathcal{V}_{p}(\lambda)$ for some $0<\lambda\leq 1$ and has the expansion of the form $(\ref{fp4eq7})$.
Then the bound
$$
|a_n(f)|\leq \frac{1-\lambda^n p^{2n}}{p^{n-1}(1-\lambda p^2)}
$$
is sharp for $n\geq 3$.
\econ

\br\label{fp4r2}
Here we note that all the results proved in \cite{BF-1} and in \cite{BF-2} for the class $\mathcal{U}_p(\lambda)$ will also
be true for the bigger function class $\mathcal{V}_{p}(\lambda)$
if we substitute $\lambda$ in place of $\lambda \mu$ and follow the same method of proof.
We also remark that the authors of \cite{OPW} has also considered similar meromorphic functions and arrive at
this conjectured bound for $|a_n|$ (compare \cite[Remark~2]{OPW}), but their study of such functions come from a different perspective.
\er

\end{document}